\documentclass[11pt,a4paper]{article}
\pdfoutput=1 
\usepackage[margin=1in]{geometry}
\usepackage{amsmath,amsthm,amssymb,enumerate, bbm ,graphicx,color,caption,upgreek, float, tikz, wasysym, subcaption,booktabs,longtable, appendix,graphics, pdfpages,rotating,mathtools,textcomp, array, blkarray}

\DeclarePairedDelimiter\floor{\lfloor}{\rfloor}
\usepackage[pdftex]{hyperref}

\definecolor{blauw}{RGB}{61,158,255}
\definecolor{donkerblauw}{RGB}{0,0,255}
\definecolor{donkergroen}{RGB}{46,148,0}
\definecolor{donkerrood}{RGB}{204,0,0}

\makeatletter 
\newcommand\mynobreakpar{\par\nobreak\@afterheading} 
\makeatother

\makeatletter
\let\@fnsymbol\@arabic
\makeatother

\newcommand{\N}{\mathbb{N}}
\newcommand{\Z}{\mathbb{Z}}

\newcommand{\R}{\mathbb{R}}

\newcommand{\F}{\mathbb{F}}

\usepackage[english]{babel}

\newtheorem{theorem}{Theorem}[section]

\newtheorem{claim}[theorem]{Claim}
\newtheorem{proposition}[theorem]{Proposition}

\theoremstyle{definition}
\newtheorem{defn}{Definition}[section]

\newtheorem*{examp*}{Example}

\DeclareMathOperator*{\wt}{wt}

\theoremstyle{plain}
\newfloat{Algorithm}{!hbt}{alg}

\newcounter{thm}[section]

\title{Uniqueness of codes using semidefinite programming} \date{}
\author{Andries E.\ Brouwer\thanks{E-mail: \href{mailto:aeb@cwi.nl}{aeb@cwi.nl}.} \, and Sven C.\ Polak\thanks{Korteweg-De Vries Institute for Mathematics, University of Amsterdam. E-mail: \href{mailto:s.c.polak@uva.nl}{s.c.polak@uva.nl}. The research leading to these
results has received funding from the European Research Council under the European Union’s Seventh Framework Programme (FP7/2007-2013) / ERC grant agreement \textnumero 339109.}}
\begin{document}
\maketitle
\setcounter{footnote}{1}

\noindent \textbf{Abstract.}  For~$n,d,w \in \N$, let~$A(n,d,w)$ denote the maximum size of a binary code of word length~$n$, minimum distance~$d$ and constant weight~$w$. Schrijver recently showed using semidefinite programming that $A(23,8,11)=1288$, and the second author that~$A(22,8,11)=672$ and~$A(22,8,10)=616$. Here we show uniqueness of the codes achieving these bounds.

Let~$A(n,d)$ denote the maximum size of a binary code of word length~$n$ and minimum distance~$d$. Gijswijt, Mittelmann and Schrijver showed that~$A(20,8)=256$. We show that there are several nonisomorphic codes achieving this bound, and classify all such codes with all distances divisible by 4.

\,$\phantom{0}$



\section{Introduction}
Let~$\F_2:=\{0,1\}$ denote the field of two elements and fix~$n \in \N$. A \emph{word} is an element~$v \in \F_2^n$. For two words~$u,v \in \F_2^n$, their \emph{(Hamming) distance}~$d_H(u,v)$ is the number of~$i$ with~$u_i \neq v_i$. A (binary) \emph{code} is a subset of~$\F_2^n$. For any code~$C$, the \emph{minimum distance}~$d_{\text{min}}(C)$ ($\in \R \cup \{\infty\}$) of~$C$ is the minimum distance between any pair of distinct code words in~$C$. The \emph{weight} $\text{wt}(v)$ of a word~$v \in \F_2^n$ is the number of nonzero entries of~$v$.  A (binary) \emph{constant weight code} is a binary code in which all code words have a fixed weight~$w$. Then $A(n,d,w)$ is defined as the maximum size of a binary constant weight~$w$ code of minimum distance at least~$d$. Moreover,~$A(n,d)$ is the maximum size of a binary code of minimum distance at least~$d$. A binary constant weight~$w$ code~$C\subseteq \F_2^n$ with~$d_{\text{min}}(C)\geq d$ is called an~$(n,d,w)$-code. A binary code~$C\subseteq \mathbb{F}_2^n$ and~$d_{\text{min}}(C) \geq d$ is called an~$(n,d)$-code.

Using semidefinite programming, some upper bounds on~$A(n,d,w)$ have recently been obtained that are equal to the best known lower bounds: it has been established that $A(23,8,11)=1288$ (see~\cite{schrijver}), and that $A(22,8,11)=672$ and $A(22,8,10)=616$ (see~$\cite{cw4}$). We show using the output of the corresponding semidefinite programs that the codes of maximum size are unique (up to coordinate permutations) for these~$n,d,w$.

For unrestricted (non-constant weight) binary codes, the bound~$A(n,d)=A(20,8)\leq256$ was obtained in~\cite{semidef}, implying that the quadruply shortened extended binary Golay code of size~$256$ is optimal.  The quadruply shortened extended binary Golay code is a linear~$(n,d)=(20,8)$-code of size~$256$ and has all distances divisible by~$4$. 

Up to equivalence there are unique $(24-i,8)$-codes of size~$2^{12-i}$ for $i=0,1,2,3$, namely the~$i$ times shortened extended binary Golay codes~\cite{brouwer2}. However, it turns out that the $4$~times shortened extended binary Golay code is not the only $(20, 8)$-code of size 256. We classify such codes with all distances divisible by~$4$, and find~$15$ such codes.

\section{The semidefinite programming upper bound}

Following~\cite{schrijver},~\cite{semidef} and~\cite{cw4}, we start by describing semidefinite programming upper bounds on~$A(n,d)$ and~$A(n,d,w)$.  Fix~$n,d,w\in \N$ and let~$N$ be either~$\mathbb{F}_2^n$ or the set of words in~$\mathbb{F}_2^n$ of weight~$w$. For~$k \in \Z_{\geq 0}$, let~$\mathcal{C}_k$ be the collection of codes~$C \subseteq N$ with~$|C|\leq k$.  We define 
\begin{align*} 
\mathcal{C}_k(D) := \{C \in \mathcal{C}_k \,\, | \,\, C \supseteq D, \, |D|+2|C\setminus D| \leq k  \}, \,\, \text{ for~$D \subseteq N$}.  
\end{align*} 
Note that then~$|C \cup C'|= |C| + |C'| -|C\cap C'| \leq 2|D| + |C \setminus D| + |C' \setminus D| -|D| \leq k$ for all~$C,C' \subseteq \mathcal{C}_k(D)$.   
Furthermore, for any function~$x : \mathcal{C}_k \to \R$ and~$D \in \mathcal{C}_k$ we define the~$\mathcal{C}_k(D) \times \mathcal{C}_k(D)$-matrix~$M_{k,D}(x)$ by $M_{k,D}(x)_{C,C'} : = x(C \cup C')$. 
\begin{align} \label{A4ndw}
A_k(n,d,w):=    \max \{ \sum_{v \in N} x(\{v\})\,\, &|\,\,x:\mathcal{C}_k \to \R_{\geq 0}, \,\, x(\emptyset )=1, x(S)=0 \text{ if~$d_{\text{min}}(S)<d$},\\[-1.1em]
& \,\, M_{k,D}(x) \succeq 0 \text{ for each~$D$ in~$\mathcal{C}_k$}\}. \notag 
\end{align}
(Here~$X \succeq 0$ means:~$X$ positive semidefinite.) Then~$A_k(n,d,w)$ is an upper bound on~$A(n,d,w)$. Similarly, one obtains an upper bound $A_k(n,d)$ on~$A(n,d)$ by setting~$N:=\mathbb{F}_2^n$ in~$\eqref{A4ndw}$, so that~$\mathcal{C}_k$ is the collection of unrestricted (not necessarily constant weight) codes of size at most~$k$. It can be proved that~$A_2(n,d)$ and~$A_2(n,d,w)$ are equal to the classical Delsarte linear programming bound in the Hamming and Johnson schemes respectively~\cite{delsarte}. 

 Let~$G$ be the set of distance preserving permutations of~$N$. In case of constant weight-codes,~$G=S_n$, where $S_n$ denotes the symmetric group on~$n$ elements, but if~$n=2w$ the group $G$ is twice as large, since then taking complements is also a distance preserving permutation of~$N$. In case of non-constant weight codes,~$G= S_2^n \rtimes S_n$, where~$S_2^n$ denotes the direct product of~$n$ copies of~$S_2$. Let~$\Omega_k$ be the set of $G$-orbits of non-empty codes in~$C_k$ and let~$\Omega_k^d \subseteq \Omega_k$ be those orbits that correspond to codes with minimum distance at least~$d$.  By averaging an optimum~$x$ over all~$x\circ g$ for~$g\in G$, one obtains the existence of a~$G$-invariant optimum solution to~$\eqref{A4ndw}$. Here~$ \circ$ denotes (function) composition, so~$x \circ g (C)= x(g(C))$ for~$C \in \mathcal{C}_k$.
 
 The original problem is equivalent to the much smaller problem in which the constraint is added that~$x$ is~$G$-invariant. We will write~$y_{\omega}$ for the common value of a $G$-invariant function~$x$ on codes~$C$ in orbit~$\omega$. Hence, the matrices~$M_{k,D}(x)$ become matrices~$M_{k,D}(y)$ and we have considerably reduced the number of variables in~$(\ref{A4ndw})$. Moreover, a block diagonalization~$M_{k,D}(y) \mapsto U_{k,D}^T M_{k,D}(y) U_{k,D}$ can be obtained reducing the sizes of the matrices involved to make the computations in~$(\ref{A4ndw})$ tractable (see~$\cite{schrijver, cw4}$ for the reductions, where we note that the reductions used in~$\cite{cw4}$ are obtained by an adaptation of the method of~$\cite{artikel}$).
 
 It can be seen (cf.~\cite{semidef, cw4}) that the nonnegativity condition on~$x$, and hence on~$y$, is already imposed by positive semidefiniteness of all matrices~$M_{k,D}(x)$. When solving the semidefinite program with a computer, we add the constraints~$y_{\omega} \geq 0$ seperately, by adding~$1\times 1$ blocks~$(y_{\omega})$ which are required to be positive semidefinite. This will allow us to easily determine which variables~$y_{{\omega}}$ will be necessarily zero in any optimum solution. This may yield necessary conditions on all optimal codes, which may give uniqueness or a classification of the optimal codes.

\subsection{Information about maximum size codes}
 Suppose that we have an instance of~$n,d$ or~$n,d,w$ for which~$A_k(n,d)=A(n,d)$ or~$A_k(n,d,w)=A(n,d,w)$, respectively. We want to obtain information about codes attaining these bounds from the semidefinite programming output. The semidefinite program~$\eqref{A4ndw}$ can be written as follows:
 \begin{align}\label{primal}
  A_k(n,d,w) = \max \left\{ \sum_{\omega \in \Omega_k^d  } b_{\omega} y_{\omega}\,\,\big{|}\,\, M= F_{\emptyset} - \sum_{\omega \in  \Omega_k^d  } F_{\omega} y_{\omega}  \succeq 0 \right\}.
 \end{align}
Here~$b_{\omega_0}=|N|$, where~${\omega_0} \in \Omega_k^d$ corresponds to the orbit of a code of size~$1$ in~$\mathcal{C}_k$, and~$b_{\omega} =0$ for all other~$\omega \in \Omega_k^d$. Moreover,~$M$ is a (large) block diagonal matrix that consists of blocks~$U_{k,D}^TM_{k,D}(y)U_{k,D}$ (which are reduced versions of the blocks~$M_{k,D}(y)$ that are required to be positive semidefinite in~$\eqref{A4ndw}$) and blocks~$(y_{\omega})$. The matrix~$F_{\omega}$ is a matrix of the same size as~$M$ with entries the coefficients of~$-y_{\omega}$ in the corresponding entries of~$M$. For each orbit~$\omega$, the matrix~$F_{\emptyset}$  is a matrix of the same size as~$M$ with entries the constant coefficients in the corresponding entries of~$M$. For two real-valued square matrices~$A,B$ of the same size, we write~$\langle A, B\rangle := \text{tr}(AB^T)$.  The dual program of~$\eqref{primal}$ then reads   
 \begin{align}\label{dual}
   \min \left\{ \langle F_{\emptyset}, X \rangle \,\,\big{|}\,\,  \langle F_{\omega}, X\rangle = b_{\omega}  \text{ for all~$\omega \in  \Omega_k^d $},\,~X\succeq 0 \right\}.
 \end{align}
If~$(M,y)$ is any optimum solution for~\eqref{primal} and~$X$ is an optimum solution for~$\eqref{dual}$ with the same value, then~$\langle M, X\rangle =0$. As~$M \succeq 0$ and~$X\succeq 0$, we have in particular~$y_{\omega} X_{\omega}=0$ for the separate~$1 \times 1$-blocks~$(y_{\omega})$ in~$M$ and~$(X_{\omega})$ in~$X$, where~$(X_{\omega})$ denotes the~$1\times 1$-block in~$X$ corresponding to the~$1\times 1$-block~$(y_{\omega})$ in~$M$. Thus
\begin{align*}
\text{$X_{\omega}>0$ in any optimum solution to~$\eqref{primal}$ } \,\,\,  \Longrightarrow \,\, \text{ $y_{\omega}=0$ in all optimum solutions to~$\eqref{dual}$}. 
\end{align*} 
 If~$y_{\omega}=0$ for all solutions to~$\eqref{A4ndw}$ with objective value~$A_k(n,d,w)=A(n,d,w)$, then for any code~$C$ of maximum size there is no subcode~$D\subset C$ with~$D \in \omega$. (Suppose otherwise, then one constructs a feasible solution to~$\eqref{A4ndw}$ by putting~$x(S)=1$ for~$S \in C_k$ with~$S \subseteq C$ and~$x(S)=0$ else, and hence by averaging over~$G$ there exists a feasible $G$-invariant solution with~$y_{\omega} >0$, a contradiction.)  So orbit~$\omega$ does not appear in any code of maximum size. Hence we know which orbits~$\omega \in \Omega_k^d$ cannot occur in a code of maximum size. We will call these orbits~\emph{forbidden} orbits. 
 
 We used the solver SDPA-GMP~$\cite{sdpa, nakata}$ to conclude which orbits are forbidden. The semidefinite programming solver does not produce exact solutions, but approximations up to a certain precision. In our case the approximations are precise enough to verify (with certainty) that certain orbits are forbidden. See the appendix for details.

\section{Self-orthogonal codes} 
If~$u,v \in \F_2^n$, we define~$(u\cap v) \in \F_2^n$ to be the word that has~$1$ at position~$i$ if and only if~$u_i=v_i=1$. The following equality is well-known and will be used often throughout the paper:
\begin{align} \label{wellknown}
d_H(u,v)= \wt(u) + \wt(v) -2 \wt(u \cap v), \,\,\, \text{ for all~$u,v \in \F_2^n$}. 
\end{align}
The function~$(u,v) \mapsto \wt(u \cap v) \pmod{2}$ is a non-degenerate symmetric $\F_2$-bilinear form on~$\F_2^n$. 
If~$(u,v)=0$, then~$u$ and~$v$ are called \emph{orthogonal}. A code~$C$ is \emph{self-orthogonal} if~$(u,v)=0$ for all~$u,v \in C$.
Given a code~$C$, the \emph{dual code}~$C^{\perp}$ is the set of all~$v \in \F_2^n$ that are orthogonal to all~$u\in C$. A code~$C$ is called \emph{self-dual} if~$C=C^{\perp}$. For small~$n$, self-dual codes are classified by Pless and Sloane~\cite{sloanepless}. 
  
 \section{Constant weight codes}
 
 With semidefinite programming three exact values of~$A(n,d,w)$ have been obtained. In~\cite{schrijver}, it is found that~$A_3(23,8,11)=1288$, matching the known lower bound and thereby proving that~$A(23,8,11)=1288$.  Similarly, in~$\cite{cw4}$, the upper bounds~$A(22,8,10)\leq 616$ and $A(22,8,11)\leq 672$ are obtained, which imply~$A(22,8,10)=616$ and~$A(22,8,11)=672$. The latter two upper bounds are in fact instances of the bound~$B_4(n,d,w)$, which is a bound in between~$A_3(n,d,w)$ and~$A_4(n,d,w)$. 
 \begin{defn}[$B_4(n,d,w)$]
 The bound~$B_4(n,d,w)$ is defined by replacing in the definition of~$A_4(n,d,w)$ from~$\eqref{A4ndw}$ the matrix~$M_{4,\emptyset}(x)$ by (the much smaller matrix)~$M_{2,\emptyset}(x)$. See~$\cite{cw4}$ for details.
 \end{defn} 
 \noindent In this section we show that the codes attaining these bounds are unique up to coordinate permutations, using the information about forbidden orbits obtained from the semidefinite programming output. In order to prove uniqueness of the~$(23,8,11)$-code of maximum size, we start by proving uniqueness of the~$(24,8,12)$-code of maximum size. The uniqueness of this code can already be obtained from the classical linear programming bound.
Below, and also later, we will need the following definition. The \emph{distance distribution}~$(a_i)_{i=0}^n$ of a code~$C \subseteq \F_2^n$ is the sequence given by~$a_i := |C|^{-1} \cdot |\{(u,v) \in C \times C \,\, | \, \, d_H(u,v)=i\}|$, for~$i=0,\ldots,n$.  The computational results we used to conclude uniqueness of the mentioned codes, are stated seperately (in~$(\ref{lpfact})$,~$(\ref{sdpfact1})$,~$(\ref{sdpfact2})$ and $(\ref{sdpfact3})$ below) at the beginning of each proof.

\subsection{\texorpdfstring{$A(24,8,12)$}{A(24,8,12)}}
\begin{proposition} \label{prop24}
Up to coordinate permutations there is a unique $(24,8,12)$-code of size~$2576$. An example is given by the set of words of weight~12 in the extended binary Golay code.
\end{proposition}
\proof 
Let~$C$ be a~$(24,8,12)$-code of size~$2576$. The classical linear programming bound in the Johnson scheme (which is equal to~$A_2(n,d,w)$) gives maximum 2576. Moreover, one has
\begin{align} \label{lpfact}
a_i =0 \text{ for }~i \notin \{0,8,12,16,24\}.
\end{align}
This information can be obtained immediately from the dual solution: the linear program contains constraints~$a_i \geq 0$. If the corresponding dual variable is~$>0$, then~$a_i=0$ in all optimum solutions to the linear program.\footnote{We used SDPA-GMP to solve this LP. The approximate dual solution allows us, with a computation similar to the computation in the Appendix in~$(\ref{error})$-$(\ref{upperbound})$ below, to give a very small $\varepsilon >0$ such that in any optimum solution,~$a_i < \varepsilon$ for~$i \notin \{0,8,12,16,24\}$. But if~$a_i>0$ for some~$(n,d,w)=(24,8,12)$-code~$C$ with~$|C|=2576$, then for this code~$a_i \geq 2 / 2576$, by definition of the distance distribution of~$C$. Since~$\varepsilon <10^{-90} < 2/2576$ for~$i \notin \{0,8,12,16,24\}$ it follows that~$a_i=0$ for any~$(24,8,12)$-code~$C$ of size~$2576$.}
 
Consider the $\F_2$-linear span~$F:=\langle C \rangle$  of~$C$. Note that~$C$, hence~$F$, is self-orthogonal, so~$|F|\leq 2^{24/2}=2^{12}$. Since~$|F|\geq |C| >2^{11}$ and~$F$ is linear, we must have~$|F|=2^{12}$, so~$F$ is self-dual. Let $u \in F$, $u \ne 0$. The sets~$\{ u + x \, |\, x \in C\} \subseteq F$ and~$C\subseteq F$ have non-empty intersection, because both sets have size~$2576> |F|/2$. So~$u+x=y$ for some~$x,y \in C$. But then~$\text{wt}(u) = d_H(x,y) \ge 8$, as~$C$ has minimum distance~$8$. It follows that~$F$ has minimum distance 8,
and we conclude that~$F$ is the extended binary Golay code. So~$C$ is the set of weight~12 words in the extended binary Golay code.
\endproof

\subsection{\texorpdfstring{$A(23,8,11)$}{A(23,8,11)}}
\begin{proposition}
Up to coordinate permutations there is a unique $(23,8,11)$-code of size~$1288$. An example is given by the set of words of weight~11 in the binary Golay code.
\end{proposition}
\proof 
Let~$C$ be a~$(23,8,11)$-code of size~$1288$. With the solution of the semidefinite program $A_3(23,8,11)$ (which is~$1288$) from~$\cite{schrijver}$ one obtains, by considering the forbidden orbits from the semidefinite programming output:\footnote{Note that the LP does not give this information: the Delsarte bound is 1417, which is not optimal.} 
\begin{align} \label{sdpfact1}
\text{if~$x,y \in C$ then~$ d_H(x,y) \leq 16$.}
\end{align}
 Construct a code~$D$ of length~$24$, weight~$12$ and size~$2576$ as follows: add a symbol~$1$ to every codeword of~$C$, put it in~$D$ and put also the complement of the resulting word in~$D$. Then~$D$ has minimum distance~$8$ by~$\eqref{sdpfact1}$. Hence~$D$ is the set of weight~$12$ words of the extended Golay code~$F$ by Proposition~$\ref{prop24}$. The automorphism group of the extended binary Golay code acts transitively on the coordinate positions~\cite{sloane}. Hence,~$C$ is the set of weight~11 words in the binary Golay code.  
\endproof

\subsection{\texorpdfstring{$A(22,8,11)$}{A(22,8,11)}}
\begin{proposition}
Up to coordinate permutations there is a unique $(22,8,11)$-code of size~$672$.
\end{proposition}
\proof 
Let~$C$ be a~$(22,8,11)$-code of size~$672$. First one concludes that~$a_{14}=0$ using the semidefinite program $B_4(22,8,11)$ from~$\cite{cw4}$. This is explained in more detail in the appendix: if~$a_{14}>0$, then~$a_{14} \geq 2/672$ and~$a_{10}+a_{14}+a_{18}+a_{22} \geq 318/672$ (by Proposition~\ref{extraprop2} below). We add these two constraints to the program~$B_4(22,8,11)$. The resulting bound is strictly smaller than~$672$, so~$a_{14}=0$ in any~$(22,8,11)$-code of size~$672$. 

Subsequently, by considering the forbidden orbits in the solution of the semidefinite program $A_3(22,8,11)$ from~$\cite{schrijver}$ with the added constraint that~$a_{14}=0$ (the solution of~$A_3(22,8,11)$ with this added constraint is~$672$) one obtains:
  \begin{align} \label{sdpfact2}
\text{if~$x,y,z \in C$ then~$d_H(x,y) \in \{0,8,12,16\}$ and~$\text{wt}(x+y+z)\in \{7,11,15\}$}.
\end{align} 
Let~$D$ be the collection of $672+672 = 1344$ codewords of length~$24$ of the form~$10x$ with~$x\in C$ together with their complements, and let $F = \langle D \rangle $ be the~$\mathbb{F}_2$-linear span of~$D$. All distances in~$D$ belong to~$\{8,12,16\}$ by~(\ref{sdpfact2}), so~$D$, and hence also~$F$, is self-orthogonal, which implies~$|F|\leq 2^{24/2} = 2^{12}$. Since all words in~$D$ have weight divisible by~$4$ and~$F$ is self-orthogonal, all words in~$F$ also have weight divisible by 4.

The code $F$ contains words of forms~$01x$,~$10y$,~$11z$ and~$00u$. Each form occurs at least~$672$ times, so $|F| \geq 4 \cdot 672 > 2^{11}$, hence $|F| = 2^{12}$ and~$F$ is self-dual.

To show that~$F$ is the extended binary Golay code, it suffices to prove that all words in~$F$ have weight~$\geq 8$, i.e., that no word in~$F$ has weight~$4$. Words of~$F$ are sums of words~$10x$ with~$x \in C$, possibly together with the all-ones word. So we must prove that sums of words~$10x$ do not have weight~$4$ or~$20$. A sum of words~$10x$ starts with~$00$ or~$10$ and is the sum of an even or odd number of words~$10x$, respectively. 

Words in~$F$ starting with~$00$ form a subcode~$F_{00}$ of~$F$ of size~$2^{10} = 1024$. If~$u \in F_{00}$, then~$\{u+10y \,|\, y \in C\} \cap \{10x  \,|\, x \in C\} \neq \emptyset$, as there are~$1024$ words in~$F$ starting with~$10$ but~$|C|=672> 1024/2$. So~$u=10x+10y$ for some~$x,y \in C$, hence~$F_{00} = \{10x + 10y \, | \, x,y \in C\}$. However, distances~$4$ and~$20$ do not occur in~$C$, so words in~$F_{00}$ do not have weight~$4$ or~$20$. 

The~$1024$ words in~$F$ starting with~$10$ are formed by the coset~$10x+ F_{00}$ (with~$x \in C$ arbitrary but fixed) and hence are a sum of three elements of the form~$10x$ with $x \in C$. But such a sum has weight~$8,12$ or~$16$ by~(\ref{sdpfact2}), implying that words in~$F$ starting with~$10$ do not have weight~$4$ or~$20$. 

Therefore weights~$4$ and~$20$ do not appear in~$F$, so~$F$ is indeed the extended binary Golay code. As the automorphism group of the extended binary Golay code~$F$ acts 2-transitively on the coordinate positions~\cite{sloane}, this implies that~$C$ is unique. 
\endproof

\subsection{\texorpdfstring{$A(22,8,10)$}{A(22,8,10)}}

\begin{proposition}
Up to coordinate permutations there is a unique $(22,8,10)$-code of size~$616$.
\end{proposition}
\proof 
Let~$C$ be a~$(22,8,10)$-code of size~$616$. First one concludes that~$a_{14}=0$ using the semidefinite program $B_4(22,8,10)$ from~$\cite{cw4}$. This is explained in more detail in the appendix: if~$a_{14}>0$, then~$a_{14} \geq 2/616$ and~$a_{10}+a_{14}+a_{18} \geq 208/616$ (by Proposition~\ref{extraprop1} below). We add these two constraints to the program~$B_4(22,8,10)$. The resulting bound is strictly smaller than~$616$, so~$a_{14}=0$ in any~$(22,8,10)$-code of size~$616$. 

Subsequently, by considering the forbidden orbits in the solution of the semidefinite program $A_3(22,8,10)$ from~$\cite{schrijver}$ with the added constraint that~$a_{14}=0$ (the solution of~$A_3(22,8,10)$ with this added constraint is~$616$) one obtains:
\begin{align} \label{sdpfact3}
\text{if~$x,y,z \in C$ then~$d_H(x,y) \in \{0,8,12,16\}$ and~$\text{wt}(x+y+z)\in \{6,10,14,22\}$}.
\end{align} 
Let~$F=\langle C \rangle$. Since~$C$ is self-orthogonal and has words of weights divisible by~$2$ but not by~$4$,~$F$ is self-orthogonal and has half of the weights divisible by~$4$ and half of the weights divisible by~$2$ but not by~$4$. Both halves of~$F$ have size~$\geq |C|=616$, but~$F$ has size~$\leq 2048$ as it is self-orthogonal. So~$|F|=2048$ and~$F$ is self-dual.

Let~$E\subseteq F$ be the subcode of~$F$ consisting of all words with weight divisible by~$4$. For each~$u \in E$, we have~$C \cap \{u+y \,|\, y \in C\} \neq \emptyset$  (as $|C|=616 >1024/2$), so~$u=x+y$ for some~$x,y \in C$. Hence~$E=\{x+y \, | \, x,y \in C\}$. By~(\ref{sdpfact3}), no word in~$E$ has weight~$4$. So weight~$4$ does not occur in~$F$.

If any word~$u$ in~$F$ has weight~$2$ then it is in~$F \setminus E = x+ E$ (with~$x \in C$ arbitrary). So it is the sum of three words in~$C$.  But such sums do not have weight~$2$ by~$(\ref{sdpfact3})$,  hence no word in~$F$ has weight~$2$. So~$F$ is a self-dual code of minimum distance~$6$.  As the self-dual~$(n,d)=(22,6)$-code is unique (cf.~$\cite{sloanepless}$),~$F$ is unique. Hence also~$C$ is unique, as it is the collection of weight~$10$ words of~$F$. (Note that two weight~$10$ words in~$F$ have distance~$0 \pmod{4}$, so distance at least~$8$, since~$d_H(u,v)=\wt(u)+\wt(v)-2\wt(u \cap v)$ for any two words~$u,v \in F$.)
\endproof

\section{Unrestricted~\texorpdfstring{$(20,8)$}{(20,8)}-codes of size 256}

Recently, Gijswijt, Mittelmann and Schrijver \cite{semidef} proved that~$A(20,8)=256$, with the semidefinite program~$A_4(n,d)$ from~$\eqref{A4ndw}$. An example of a code attaining this bound is the four times shortened extended binary Golay code, which has distance distribution
\begin{align}\label{distgol}
a_0=1, \,\, a_8=130,\,\, a_{12}=120, \,\,a_{16}=5,\,\,\text{$a_i=0$ for all other~$i$}. 
\end{align}
This code is formed by the words starting with~$0000$ in the extended binary Golay code with these first four coordinate positions removed. 

Two binary codes~$C,D \subseteq \F_2^n$ are \emph{equivalent} if~$D$ can be obtained from~$C$ by first permuting the~$n$ coordinates and by subsequently permuting the alphabet~$\{0,1\}$ in each coordinate separately.

Up to equivalence there are unique $(24-i,8)$-codes of size~$2^{12-i}$ for $i=0,1,2,3$, namely the~$i$ times shortened extended binary Golay codes \cite{brouwer2}. In this section we show that there exist several nonisomorphic~$(20,8)$-codes of size~$256$. First we show that there exist such codes with different distance distributions. Subsequently we classify such codes with all distances divisible by~$4$.

We start by recovering information about possible distance distributions from the semidefinite program~$A_4(20,8)$. Write~$\omega_t \in \Omega_k$ for the orbit of two words at Hamming distance~$t$.  From a code~$C$ with distance distribution~$(a_i)$, one constructs a feasible solution to~$\eqref{A4ndw}$ by putting~$x(S)=1$ for~$S \in C_k$ with~$S \subseteq C$ and~$x(S)=0$ else, and hence by averaging over~$G$ one obtains a feasible $G$-invariant solution with variables~$y_{\omega}$. This solution has
\begin{align*}
y_{\omega_t} &= \frac{1}{|G|} \sum_{g \in G} x\circ g(\{x,y\}) =  \frac{t!(20-t)!}{|G|}|\{(u,v) \in C^2 \,\, : \,\, d_H(u,v)=t \}|  \notag 
\\&= \frac{|\{(u,v) \in C^2 \,\, : \,\, d_H(u,v)=t \}|}{2^{20} \binom{20}{t}}   = \frac{|C| a_t}{2^{20} \binom{20}{t}}, 
\end{align*}
where~$\{x,y\}$ is any pair of words with distance~$t$ and~$G=S_2^{20} \rtimes S_{20}$. So we can add linear constraints on the~$a_i$ as linear constraints on the variables~$y_{\omega_t}$ to our semidefinite program.

The inner distribution~$(a_i)$ is not determined uniquely by the requirement that it is an optimal solution of
the semidefinite program~$A_4(20,8)$ from~$\eqref{A4ndw}$.\footnote{By contrast, in all constant weight cases considered in this paper, the values of the~$y_{\omega_t}$ give the unique distance distribution of the (unique up to coordinate permutations)~$(n,d,w)$-codes of maximum size.} We find minimum possible values for some of the~$a_i$ for the case where all distances are even as follows. For any code~$C$, the~$a_i$ ($i\neq 0$) are integer multiples of~$2/|C|$. So for any~$(20,8)$-code of size~$256$, if~$a_{16}<1$ then~$a_{16}\leq 254/256$.  With the constraint~$a_{16} \leq 254/256$ the semidefinite program returns an objective value strictly smaller than~$256$. So~$a_{16} \geq 1$. Similarly, we find~$a_8 \geq 126$ and~$a_{12} \geq 96$. If we simultaneously add the constraints~$a_8 \leq 126$, $a_{12} \leq 96$ and~$a_{16} \leq 1$, the semidefinite program returns~$256$ as objective value, and the values of~$y_{\omega_{10}}$ and~$y_{\omega_{14}}$ force~$a_{10}=a_{14}=16$. Therefore, apart from the 4 times shortened extended binary Golay code, also a code with
\begin{align}\label{dist20}
a_1=1,\,\,a_{8}=126,\,\, a_{10}=16, \,\, a_{12}=96,\,\, a_{14}=16,\,\, a_{16}=1,\,\,\text{$a_i=0$ for all other~$i$},
\end{align}
is allowed by the program~$A_4(20,8)$. Such a code exists, as the following construction demonstrates. 

Start with the extended binary Golay code~$F$ containing the weight~$8$ word~$u$ with all 1s in the first eight positions. As~$A(24-8,8)=A(16,8)=32$ (see~$\cite{brouwertableand}$), there can be at most~$32$ words in~$F$ starting with~$8$ zeroes. These form a linear subcode~$E$ of~$F$. As any word in~$F$ has an even number of~$1$s at the first eight positions, there are at most~$2^7=128$ distinct cosets~$E+v$ in~$F$. As~$32 \cdot 128=2^{12}=|F|$ it follows that~$|E|=32$ and there are exactly~$128$ distinct cosets~$E+v$.
 
 So if we specify a string of~$8$ symbols with an even number of ones and take all words in~$F$ having these~$8$ fixed symbols in the first~$8$ positions, we obtain a subcode~$D$ of~$F$ of size~$32$ and minimum distance at least~$8$.  Choose the following~$8$ specifications, each giving a subcode~$D$ of size~$32$ and minimum distance~$8$. 
\begin{align*}
    \begin{matrix} 
00000000 \\
11000000 \\
10100000 \\
10010000 \\
10001000 \\
10000100 \\
10000010 \\
10000001
\end{matrix} \,\,\, \quad \text{ and then replace the first~$8$ coordinates by }\quad \,\,\begin{matrix} 
0000 \\
1100 \\
1010 \\
1001 \\
0110 \\
0101 \\
0011 \\
1111
\end{matrix}.
\end{align*}
This yields a~$(20,8)$-code of size~$8\cdot 32= 256$ in which distances~$10$ and~$14$ occur. Note that this code indeed has minimum distance at least~$8$: first observe that each code~$D$ has minimum distance at least~$8$. Then note that for two different specifications the first part (the first~$8$ positions) had distance at most~$2$ before the replacement, so the second part has distance at least~$6$. After the replacement of the first part, the first~$4$ positions have distance at least~$2$, so two words obtained from different specifications have, after the replacement, in total distance at least~$2+6=8$.  One verifies by computer that its distance distribution is given by~$\eqref{dist20}$. So, there exist~$(20,8)$-codes of maximum size with distance distribution~$\eqref{distgol}$ as well as with~$\eqref{dist20}$. 

\subsection{Unrestricted~\texorpdfstring{$(20,8)$}{(20,8)}-codes of maximum size with all distances divisible by~4}

In this section we give a classification of the~$(20,8)$-codes of size~$256$ with all distances divisible by~$4$. An example of such a code is the quadruply shortened extended binary Golay code~$B$, which is linear. There is, up to equivalence, only one such code, since the automorphism group of the extended Golay code acts 5-transitively on the coordinate positions~\cite{brouwer2}. (Moreover, Dodunekov and Encheva~$\cite{dodunekov}$ have proved that there exists, up to equivalence, only one \emph{linear}~$(20,8)$-code of size~$256$.) The shortened extended binary Golay code contains~$5$ words of weight~$16$, forming a subcode~$D$. Since the minimum distance is~$8$, each of those words must have the~$0$'s at different positions. So we can assume that 
\begin{align} \label{Dmat}
D=     \begin{matrix} 
00001111111111111111 \\
11110000111111111111 \\
11111111000011111111 \\
11111111111100001111 \\
11111111111111110000
\end{matrix},
\end{align}
 with linear span~$\langle D \rangle \subseteq B$ of dimension~$4$. So~$B$ is a union of~$16$ cosets~$u+\langle D\rangle $. If we replace a coset~$u + \langle D \rangle$ by its complement~$\mathbf{1}+u + \langle D \rangle$, we obtain another~$(20,8)$-code of maximum size that is not linear, so this is really a different code. Note that all distances remain divisible by --but not equal to-- $4$, as~$d_H(x,y)\in \{8,12\}$ for any~$x \in u + \langle D \rangle$ and~$y \in B \setminus (u + \langle D \rangle) $, so~$d_H(\mathbf{1}+x,y) \in \{8,12\}$. By replacing any of the~$16$ cosets~$u+\langle D \rangle$ in~$B$ with~$\mathbf{1}+u+\langle D \rangle$, we obtain~$2^{16}=65536$ codes with all distances divisible by~$4$. 
 
 In this section we will first prove that any maximum-size~$(20,8)$-code  with all distances divisible by~$4$ is equivalent to one of the~$2^{16}$ thus obtained codes. Secondly, we will obtain (by computer) that these~$2^{16}$ codes can be partitioned into~$15$ equivalence classes. 
 
 In order to prove the first result, we start by proving two auxiliary propositions. Let~$C$ be any~$(20,8)$-code of size~$256$ with all distances divisible by~$4$ and containing~$\mathbf{0}$, the zero word. Define~$E:=\langle C, \mathbf{1} \rangle$ to be the linear span of~$C$ together with the all-ones vector.
\begin{proposition} \label{Eprop}
Up to a permutation of the coordinate positions, the codes~$E$ and $\langle B, \mathbf{1}\rangle $ are the same.
\end{proposition}
\proof 
 After the constraints~$a_i=0$ if~$4 \nmid i$ and~$a_{20} \geq 2/256$ are added to the ordinary LP-bound for~$(n,d)=(20,8)$, the linear program returns a solution strictly smaller than~$256$. Therefore~$a_{20}=0$ in any~$(20,8)$-code of size~$256$ with all distances divisible by~$4$. As~$A(20,8,8)=130$ (cf.~\cite{brouwertable}) and $A(20,8,4)=\floor{20/4}=5$, one has~$a_{8}\leq 130$ and~$a_{16} \leq 5$. Moreover, the LP-bound contains the inequalities $a_8-a_{12}-3a_{16}+5 \geq 0$ and~$-a_{8} -a_{12} +31 a_{16} +95 \geq 0$  (given that~$a_{20}=0$). We add those two equations and use that~$a_{16} 
 \leq 5$ to obtain that~$a_{12}\leq 120$. As~$256=1+a_8+a_{12}+a_{16}$, the distance distribution of any~$(20,8)$-code of size 256 with all distances divisible by 4 is given by~$(\ref{distgol})$.  Moreover, the values in~$(\ref{distgol})$ are not mere averages: as~$A(20,8,8)=130$ and $A(20,8,4)=5$, the number of words at distance~$i$ from \emph{any} word~$u \in C$ is specified by~$(\ref{distgol})$.

Since~$C$ is self-orthogonal and has all distances divisible by~$4$, also~$E$ is self-orthogonal and has all distances divisible by~$4$. Furthermore,~$E$ has dimension~$9$. To see this, note that~$\mathbf{1} \notin C$ since~$a_{20} = 0$, so ~$|E|\geq 257$, so~$|E| \geq 512$. On the other hand,~$\dim E <10$, as~$E$ is self-orthogonal with all distances divisible by~$4$, but there does not exist a self-dual code of length~$20$ with all distances divisible by~$4$ (cf.~$\cite{pless}$). So~$\dim E =9$ and~$|E|=512$, implying that
\begin{align} \label{ofof}
\text{for every word~$u \in E$ one has~$u\in C$ or~$\mathbf{1}+u \in C$.} 
\end{align}
For any code we write~$A_i$ for the number of words of weight~$i$. Since~$C$ has weights $A_0=1$, $A_8=130$, $A_{12}=120$, $A_{16}=5$, we conclude that~$E$ has weights 
\begin{align*} 
A_0=1,\, A_4=5, \,A_8=250,\, A_{12}=250, \,A_{16}=5,\, A_{20}=1.
\end{align*}
The orthogonal complement~$E^\perp$ of~$E$ has dimension 11, and is a union $E \cup (a+E) \cup (b+E) \cup (c+E)$.
Here $a,b,c$ have even weight (because $\mathbf{1} \in E$), so each of $\langle a,E\rangle $ and $\langle b,E\rangle $
and $\langle c,E \rangle$ is self-dual. This means that $a,b,c$ are mutually non-orthogonal.

Look at Pless~$\cite{pless}$ to find the self-dual codes of length~$n=20$ and dimension~$10$.
There are 16 such codes, but we can forget about those with $A_8 < 250$.
There is a unique self-dual code of length~$n=20$ and dimension~$10$ with $A_8 \geq 250$, namely $M_{20}$ with weight enumerator 
\begin{align*} 
A_0 = 1,\, A_4 = 5,\, A_6 = 80,\, A_8 = 250,\, A_{10} = 352, \,A_{12}=250,\, A_{14}=80,\, A_{16}=5,\, A_{20}=1,
\end{align*} 
and~$E$ is the subcode of $M_{20}$ consisting of the words of weight divisible by~$4$, hence is unique. This means that~$E$ does not depend on the particular choice of~$C$ (up to a permutation of the coordinates), proving the desired result. 
\endproof 

\begin{proposition} \label{Dprop}
Let~$C$ be any~$(20,8)$-code of size~$256$ with all distances divisible by~$4$ containing~$\mathbf{0}$. Then~$C$ is invariant under translations by weight~$16$ words from~$C$.
\end{proposition}
\proof 
Clearly, if $a,b,c \in \F_2^{20}$ with $\wt(a)=16$ and $a+b+c=\mathbf{1}$, then $d_H(b,c) = 4$.
Now let~$b \in C$ be arbitrary, and~$a \in C$ a weight~$16$ vector. Then we have~$a+b+c \neq \mathbf{1}$ for all~$c \in C$. So~$\mathbf{1}+a+b \notin C$ while~$\mathbf{1}+a+b \in \langle C,\mathbf{1} \rangle$. Hence,~$a+b \in C$, by~$\eqref{ofof}$. So indeed,~$C$ is invariant under translations by weight~$16$ words from~$C$. (This in particular implies that~$C$ contains the $4$-dimensional linear span of the weight~$16$ vectors.)
\endproof 
Fix representatives~$u_1,\ldots,u_{16}$ for the cosets~$u_i+\langle D \rangle$ of the linear quadruply shortened binary extended Golay code~$B$ (see Table~$\ref{tab1}$ for a possible choice). Using Propositions~$\ref{Eprop}$ and~$\ref{Dprop}$, we obtain the first main result of this section.

\begin{proposition}\label{216}
Let~$C$ be any~$(20,8)$-code of size~$256$ with all distances divisible by~$4$. Then~$C$ is equivalent to~$B$ with some of the cosets~$u_i+\langle D \rangle$ replaced by~$\mathbf{1}+u_i+\langle D \rangle$. 
\end{proposition}
\proof 
 By applying a distance preserving permutation to~$C$, which is possible by Proposition~$\ref{Eprop}$, we may assume that~$C$ contains~$\mathbf{0}$ and that~$\langle C,\mathbf{1} \rangle = \langle B, \mathbf{1} \rangle$. Then~$C$ contains~$5$ weight~$16$ vectors of the form~$\eqref{Dmat}$. By Proposition~$\ref{Dprop}$, the code~$C$ is a union of~$16$ cosets~$u+\langle D \rangle$, for some vectors~$u$. The code~$ \langle B, \mathbf{1} \rangle=\langle C, \mathbf{1} \rangle $ is a union of cosets~$u_i+\langle D \rangle $ together with their complements~$\mathbf{1}+u_i + \langle D \rangle$. This implies by~$\eqref{ofof}$ that each coset of~$C$ has the form~$u_i+\langle D \rangle$ or~$\mathbf{1}+u_i+\langle D \rangle$ (and~$C$ cannot contain both~$u_i$ and~$\mathbf{1}+u_i$ at the same time as~$a_{20}=0$), as required. 
\endproof

\begin{table}[ht]
{\small
  \setlength{\abovedisplayskip}{6pt}
  \setlength{\belowdisplayskip}{\abovedisplayskip}
  \setlength{\abovedisplayshortskip}{0pt}
  \setlength{\belowdisplayshortskip}{3pt}
  \setlength{\tabcolsep}{.16666667em}
\begin{align*}
{\begin{tabular}{|l| *{15}{c}|}\hline
   Coset representative &  $C_1$ & $C_2$ & $C_3$ & $C_4$ & $C_5$ & $C_6$  & $C_7$& $C_8$ & $C_9$ & $C_{10}$ & $C_{11}$ & $C_{12}$ & $C_{13}$  & $C_{14}$ & $C_{15}$    \\
\hline
$u_{1\phantom{0}} = 00000000000000000000$ \,\,\,\,\,\,\,\,\,\,    & &1&1& 1&1&1&   1&1&1&1&  & & & 1&    \\
$u_{2\phantom{0}} = 00000101010101011010$                         & & &1& 1&1&1&   1&1&1&1& 1&1&1&  &1   \\
$u_{3\phantom{0}} = 00001001011001101100$                         & & & & 1&1& &   1& & &1& 1&1&1& 1&    \\
$u_{4\phantom{0}} = 00001100001100110110$                         & & & &  &1& &    & & & & 1& & &  &1   \\
$u_{5\phantom{0}} = 10100000010101101001$                         & & & &  & &1&   1&1&1&1& 1&1&1&  &    \\
$u_{6\phantom{0}} = 10100101000000110011$                         & & & &  & & &    &1& & &  &1& & 1&1   \\
$u_{7\phantom{0}} = 10101001001100000101$                         & & & &  & & &    & &1& &  & &1& 1&1   \\
$u_{8\phantom{0}} = 10101100011001011111$                         & & & &  & & &    & & & &  & & &  &    \\
$u_{9\phantom{0}} =  11000000011000110101$                        & & & &  & & &    & & &1& 1&1&1& 1&1   \\
$u_{10} = 11000101001101101111$                                   & & & &  & & &    & & & &  & &1& 1&1   \\
$u_{11} = 11001001000001011001$                                   & & & &  & & &    & & & &  & & &  &    \\
$u_{12} = 11001100010100000011$                                   & & & &  & & &    & & & &  & & &  &    \\
$u_{13} = 01100000001101011100$                                   & & & &  & & &    & & & &  & & &  &    \\
$u_{14} = 01100101011000000110$                                   & & & &  & & &    & & & &  & & &  &    \\
$u_{15} = 01101001010100110000$                                   & & & &  & & &    & & & &  & & &  &    \\
$u_{16} = 01101100000001101010$                                   & & & &  & & &    & & & &  & & &  &    \\
\hline  
\end{tabular}               
              }
  \end{align*}
}%
\caption{\label{tab1}\small The~$(20,8)$-codes of size~$256$ with all distances divisible by~$4$. The quadruply shortened extended binary Golay code~$B=C_1$ is the union~$\cup_{i=1}^{16} (u_i+ \langle D \rangle)$. The other codes~$C_j$ ($j=2,\ldots,15$) are obtained from~$B$ by replacing the coset~$u_i+\langle D \rangle$ by~$\mathbf{1}+u_i + \langle D \rangle$ if there is a~$1$ in entry~$(u_i,C_j)$ in the above table. }
\end{table}

It remains to classify the~$2^{16}=65536$ codes obtained from~$B$ by replacing some of the cosets~$u_i+\langle D \rangle$ by~$\mathbf{1}+u_i+\langle D \rangle$. For this we use the graph isomorphism program \texttt{nauty}~$\cite{dreadnaut}$. For any code~$C$ of word length~$n$ containing~$m$ codewords, a graph with~$2n+m$ vertices is created: one vertex for each codeword~$u \in C$ and two vertices~$0_i$ and~$1_i$ for each coordinate position. Each code word~$u$ has neighbor~$0_i$ if~$u_i=0$ and~$1_i$ if~$u_i=1$ ($i=1,\ldots,n)$. Moreover, there are edges~$\{0_i,1_i\}$ ($i=1,\ldots,n$). 

All code words have degree~$n$ and the coordinate positions have (in this case) larger degree. An automorphism of this graph permutes the codewords and permutes the coordinate positions. In this way one finds a subgroup of~$S_2^n \rtimes S_n$ that fixes~$C$ and the question of code equivalence is transformed into a question of graph isomorphism. With the program \texttt{nauty} we compute a canonical representative for each of the~$2^{16}$ mentioned codes. In this way we find that the~$2^{16}$ codes from Proposition~$\ref{216}$ can be partitioned into~$15$ equivalence classes. See Table~$\ref{tab1}$ for the classification.

\begin{proposition}
There are $15$ different $(20,8)$-codes of size~$256$ with all distances divisible by~$4$ up to equivalence. \qed
\end{proposition}

 \section*{Appendix: approximate solutions}
 
 The semidefinite programming solver does not produce exact solutions, but approximations up to a certain precision. Here we show that we have enough precision to conclude~$\eqref{sdpfact1}$,~$\eqref{sdpfact2}$ and~$\eqref{sdpfact3}$. 
 
 Let~$(M,y)$ be feasible for~$\eqref{primal}$ with optimum value at least~$A(n,d,w)$ and suppose that~$X$ is an approximation of the dual program. That is, we have~$X \succeq 0$\footnote{We used a seperate java program to verify that~$X \succeq 0$ (in fact,~$X \succ 0$) in the SDP-outputs used in this paper.} and for all~$\omega \in \Omega_k^d$:
\begin{align} \label{error}
    \langle X, F_{\omega} \rangle = b_{\omega} + \epsilon_{\omega}, 
\end{align}
 for small~$\epsilon_{\omega}$.  Consider one particular~$\omega \in \Omega_k^d$. There is a~$1 \times 1$-block~$(y_{\omega})$ in~$M$ and hence also a corresponding~$1 \times 1$-block~$(X_{\omega})$ in~$X$. Remove these~$1\times1$ blocks from~$M$ and~$X$ and call the resulting matrices~$M'$ and~$X'$. Then
 \begin{align}
     0 &\leq \langle M', X' \rangle = \langle M,X\rangle -y_{\omega}X_{\omega}  \notag 
     \\&= \langle F_{\emptyset},X\rangle - \sum_{\omega \in \Omega_k^d} y_{\omega}b_{\omega}  + \sum_{\omega \in \Omega_k^d} y_{\omega} \epsilon_{\omega} -X_{\omega} y_{\omega}, \label{errorcomp}
 \end{align}
 as~$M'$ and~$X'$ are positive semidefinite, where we used~$\eqref{error}$ and the definition of~$M$ from~$\eqref{dual}$ in the second equality. Note that~$\sum_{\omega \in \Omega_k^d} y_{\omega}b_{\omega}$ is bounded from below by~$|C|$, and each~$y_{\omega}$ is bounded from above by~$y_0$ (this can be done since all~$2 \times 2$ principal submatrices in the semidefinite program~$(\ref{primal})$ are positive semidefinite), which is bounded from above by~$1$. Hence
 \begin{align} \label{errorcomp2}
 X_{\omega}y_{\omega} \leq \langle F_{\emptyset},X \rangle - |C| +  \sum_{\omega \in \Omega_k^d}\epsilon_{\omega}.
 \end{align}
 The numbers~$\epsilon_{\omega}$ are easily calculated from the dual solution, just as the dual approximate objective value~$\langle F_{\emptyset},X\rangle$. So we find  a constant~$c_{\omega}$ from the semidefinite programming dual approximation~$X$ such that
 \begin{align}
 X_{\omega}y_{\omega} \leq c_{\omega}.
 \end{align}
In the case of~$A(23,8,11)$  one can conclude with the semidefinite program~$A_3(23,8,11)$, which can be solved with SDPA-GMP~$\cite{sdpa, nakata}$ within minutes, that
\begin{align} \label{upperbound} 
y_{\omega} \leq 10^{-90}
\end{align} 
for all orbits corresponding to codes not satisfying~$(\ref{sdpfact1})$. Let~$\omega$ be an orbit for which~$(\ref{upperbound})$ holds. If there exists a code~$C$ of maximum size containing a subcode~$D\subset C$ with~$D \in \omega$, then one constructs a feasible solution to~$\eqref{A4ndw}$ by putting~$x(S)=1$ for~$S \in C_k$ with~$S \subseteq C$ and~$x(S)=0$ else, and hence by averaging over~$G$ there exists a feasible $G$-invariant solution with~$y_{\omega} \geq 1/|G|$ (this lower bound is not best possible, but sufficient). In our case,~$G=S_{23}$, so~$y_{\omega} \geq 1/23! > 10^{-23}$, which gives a contradiction with~$\eqref{upperbound}$. In this way, one verifies that all orbits not satisfying~$(\ref{sdpfact1})$ are forbidden, thereby establishing~$(\ref{sdpfact1})$. We used a seperate java program to check that~$X\succeq 0$ (in fact,~$X \succ 0$) and to compute the error terms as in~$(\ref{error})$ and~$(\ref{errorcomp2})$.

Next, we consider the cases~$A(22,8,10)=616$ and~$A(22,8,11)=672$. First we assume that~$a_{14}=0$ for all maximum-size~$(22,8,10)=616$ and~$(22,8,11)$-codes~$C$.  We write~$\omega_t$ for the orbit of two words at Hamming distance~$t$.  Adding the constraint~$y_{\omega_{14}}=0$ to the programs~$A_3(n,d,w)$  for these two cases of~$n,d,w$ gives~$A_3(n,d,w)=A(n,d,w)$. In this way one shows, in the same way as in the previous paragraph, that all orbits not satisfying~$(\ref{sdpfact2})$ and~$(\ref{sdpfact3})$ are forbidden, provided that~$a_{14}=0$. So in order to establish~$(\ref{sdpfact2})$ and~$(\ref{sdpfact3})$, it remains to prove that 
\begin{align}\label{rtp}
\text{if~$C$ is a maximum-size~$(22,8,10)$- or~$(22,8,11)$-code, then~$a_{14}=0$.} 
\end{align} 
Suppose to the contrary that~$C$ is a code as in~$\eqref{rtp}$, yet~$a_{14}>0$. Then~$a_{14} \geq 2/|C|$. We will show that this is not possible by adding constraints to the (large) program~$B_4(n,d,w)$. The semidefinite program will then give~$\floor{B_4(n,d,w)} < A(n,d,w)$ and we will arrive at a contradiction, as~$B_4(n,d,w)$ is an upper bound for~$A(n,d,w)$. To find a better lower bound on some of the~$a_i$, we use the following two propositions. We use in both propositions that for two words~$u,v$ in a constant weight~$w$ code~$C$,
\begin{align} \label{oddeven}
d_H(u,v) \equiv 2 \pmod{4} \,\,\,\,\, \Longleftrightarrow \,\,\,\,\, \wt(u \cap v) \not\equiv w \pmod{2},
\end{align} 
which follows from~\eqref{wellknown}. 
In the next two propositions we will call~$\wt(u\cap v)$ the \emph{inner product} of~$u$ and~$v$.
\begin{proposition} \label{extraprop1}
Let~$C$ be a~$(22,8,10)$-code of size~$616$ with~$a_{14} \geq 2 / 616$. Then~$a_{10}+a_{14}+a_{18} \geq 208/616$. 
\end{proposition}
\proof 
Suppose that~$a_{10}+a_{14}+a_{18} < 208/616$ (note that~$a_{22}=0$ as two weight~$10$ words cannot have distance~$22$). Then, by~$\eqref{oddeven}$, there are at most~$206/2=103$ pairs of words in~$C$ with odd inner product. Let~$\{b,b'\} \subseteq C$ be such a pair of words. Starting with~$b\in C$, and greedily picking vectors~$d \in C$ such that the inner product of~$d$ with the already chosen vectors is even, we end with a self-orthogonal subcode~$B$ of~$C$ of size~$\geq 616-103=513$. Starting with~$b'$, we repeat the same process to end up with a self-orthogonal subcode~$B'$ of~$C$ (containing~$b'$) of size~$\geq 513$. Furthermore,~$D:= B \cap B'$ has~$|D| \geq 512$. (To see this, note that $b,b'\notin D$ with odd inner product. Every word~$v \in C \setminus (D\cup \{b,b'\})$ has odd inner product with some word in~$C$ so there are at most~$103-1=102$ of such words~$v$, as~$b,b'$ is already a pair with odd inner product.) Write~$F:=\langle D \rangle$. Then~$F$ is self-orthogonal, as~$D$ is self-orthogonal.

Since~$\langle B \rangle$ and~$\langle B' \rangle$ are self-orthogonal codes, they have dimension at most~$11$. Since~$\langle B \rangle \neq \langle B'\rangle$, as~$b \notin \langle B'\rangle $, we have~$\dim F =  \dim \langle B  \cap B'\rangle \leq \dim( \langle B \rangle \cap \langle B'\rangle) \leq 10$. On the other hand, we have~$ |F| >512$ (as~$D \subseteq F$ and the zero word is contained in~$F$). So~$\dim F=10$. Moreover,~$F$ has minimum distance~$8$. To see this, let~$u$ be any nonzero word in~$F \setminus D$. Then~$D \cap (u +D) \neq \emptyset$, as both~$D$ and~$u+D$ have size~$512$ and are contained in~$F \setminus \{ \mathbf{0}\}$, a set of size~$1023$. So~$u$ is the sum of two words of~$D$, and hence has weight at least~$8$ (as~$D$ has minimum distance at least~$8$).

So~$F$ is a self-orthogonal code of word length~$22$, dimension~$10$ and minimum distance~$8$. Such a code is the twice shortened extended binary Golay code (see~$\cite{brouwer2}$), which does not contain words of weight~$10$. But all words in~$D \subseteq F$ have weight~$10$, a contradiction.
\endproof 

\begin{proposition} \label{extraprop2}
Let~$C$ be an~$(22,8,11)$-code of size~$672$ with~$a_{14} \geq 2 / 672$. Then~$a_{10}+a_{14}+a_{18}+a_{22} \geq 318/672$. 
\end{proposition}
\proof 
Suppose that~$a_{10}+a_{14}+a_{18}+a_{22} < 318/616$. Then, by~$\eqref{oddeven}$, there are at most~$316/2=158$ pairs of words in~$C$ with even inner product. Let~$\{b,b'\} \subseteq C$ be such a pair of words. Starting with~$b\in C$, and greedily picking vectors~$d \in C$ such that the inner product of~$d$ with the already chosen vectors is odd, we end with a subcode~$B$ of~$C$ of size~$\geq 672-158=514$. Now, add an extra symbol~$1$ to every codeword in~$B$, to obtain a self-orthogonal code~$D$ of length~$23$. As~$D$ is self-orthogonal,~$\dim \langle D \rangle \leq \floor{23/2}=11$. 

Starting with~$b' \in C$, we repeat the same process to end up with a subcode~$B'$ of~$C$ (containing~$b'$) of size~$\geq 514$ such that all pairs of words in~$B'$ have odd inner product. Add an extra symbol~$1$ to every code word in~$B'$ to obtain a self-orthogonal code~$D'$ of length~$23$, so~$\dim \langle D'\rangle \leq 11$. Note that~$\langle D \rangle  \neq \langle D'\rangle$, as~$1b' \notin \langle D\rangle $.

Furthermore,~$E:= D \cap D'$ has~$|E| \geq 513$ and all words start with~$1$. (To see this, note that $b,b'\notin B \cap B' $ with even inner product. Every word~$v \in C \setminus ((B \cap B' ) \cup \{b,b'\})$ has even inner product with some word in~$C$ so there are at most~$157$ of such words~$v$.)  Hence~$|\langle E\rangle | \geq 2 \cdot 513$, so~$\dim \langle E\rangle  \geq 11$. But~$\langle D\rangle $ and~$\langle D'\rangle $ are distinct codes of dimension~$\leq 11$, so their intersection has dimension~$<11$, hence
$$
\dim\langle E \rangle=\dim\langle D \cap D'\rangle \leq \dim (\langle D\rangle \cap \langle D' \rangle) <11,
$$
a contradiction. 
\endproof 

From a code~$C$ with distance distribution~$(a_i)$, one constructs a feasible solution to~$\eqref{A4ndw}$ by putting~$x(S)=1$ for~$S \in C_k$ with~$S \subseteq C$ and~$x(S)=0$ else, and hence by averaging over~$G$ there exists a feasible $G$-invariant solution. This solution has
\begin{align*}
y_{\omega_t} &= \mbox{$\frac{1}{|G|}$} \sum_{g \in G} x\circ g(\{x,y\}) =  \frac{(\frac{t}{2})!(w-\frac{t}{2})!(\frac{t}{2})!(22-w-\frac{t}{2})!}{|G|}|\{(u,v) \in C^2 \,\, : \,\, d_H(u,v)=t \}|  \notag 
\\&= \frac{|\{(u,v) \in C^2 \,\, : \,\, d_H(u,v)=t \}|}{\binom{22}{w} \cdot\binom{22-w}{t/2} \binom{w}{t/2}}   = \frac{|C| a_t}{\binom{22}{w} \binom{22-w}{t/2} \binom{w}{t/2}} = \frac{y_{\omega_0} a_t}{\binom{22-w}{t/2} \binom{w}{t/2}}, 
\end{align*}
where~$\{x,y\}$ is any pair of constant-weight~$w$ words with distance~$t$ and~$G=S_{22}$, the symmetric group on 22 elements. 

So we can add linear constraints on the~$a_i$ as linear constraints on the variables~$y_{\omega_t}$ to our semidefinite program. To the program~$B_4(22,8,10)$ we add the constraints~$a_{14} \geq 2/616$ and~$a_{10}+a_{14}+a_{18} \geq 208/616$. To~$B_4(22,8,11)$ we add the constraints~$a_{14} \geq 2/672$ and~$a_{10}+a_{14}+a_{18}+a_{22} \geq 318/672$. Write~$B_4^*(n,d,w)$ for the resulting bound after adding these constraints. We find~$B_4^*(n,d,w)<A(n,d,w)$ in both cases (which we verified using the dual solution), which is not possible.\footnote{The SDP-solutions show $B_4^*(22,8,11)< 671.885 <672$ and~$B_4^*(22,8,10)<615.935<616$.} This establishes~$(\ref{rtp})$ and hence completes the verification of~$(\ref{sdpfact2})$ and~$(\ref{sdpfact3})$.

The time needed to solve the semidefinite programs~$B_4^*(22,8,11)$ and~$B_4^*(22,8,10)$ varied from one to three weeks with sufficient precision to conclude that~$B_4^*(n,d,w) < A(n,d,w)$ in these two cases (with SDPA-DD). By contrast, the semidefinite programs for~$A_3(22,8,10)$ and~$A_3(22,8,11)$ can be solved withvery high precision within minutes (with SDPA-GMP). 

The computer programs we used to generate input for the SDP-solver can be found in~\cite{programs}. Also, the input and output files for the SDP solver can be found in this folder, and a java program to inspect the outputs.

\section*{Acknowledgements}
The second author wants to thank Lex Schrijver for useful discussions. Furthermore, we would like to thank the referees for their valuable comments to improve the paper.

\selectlanguage{english}


\begin{thebibliography}{99}
  \addcontentsline{toc}{chapter}{References}
\bibliographystyle{alpha}
\small 
\bibitem{brouwertable} A.\ E.\ Brouwer, Tables with bounds on~$A(n,d,w)$, \url{http://www.win.tue.nl/~aeb/codes/Andw.html}.

\bibitem{brouwertableand} A.\ E.\ Brouwer, Tables with bounds on~$A(n,d)$, \url{http://www.win.tue.nl/~aeb/codes/binary-1.html}.

\bibitem{brouwer2} A.\ E.\ Brouwer, Block designs, in: \textsl{Handbook of Combinatorics},
R.\ Graham, M.\ Groetschel, L.\ Lov\'asz, eds.,
Elsevier, 1995, 693--745.

\bibitem{delsarte} P.\ Delsarte, An algebraic approach to the association schemes of coding theory, \textsl{Philips Res.\ Repts.\ Suppl.} No.\ 10, 1973.

\bibitem{dodunekov} S.\ M.\ Dodunekov, S.\ B.\ Encheva, Uniqueness of Some Linear Subcodes of the Extended Binary Golay Code, \textsl{Probl.\ Peredachi
Inf.},  29 (1993), 45-–51

\bibitem{semidef}D.\ C.\ Gijswijt, H.\ D.\ Mittelmann, A.\ Schrijver, Semidefinite code bounds based on quadruple distances, \textsl{IEEE Transactions on Information Theory},  58 (2012), 2697--2705.

\bibitem{artikel}B.\ M.\ Litjens, S.\ C.\ Polak, A.\ Schrijver, Semidefinite bounds for nonbinary codes based on quadruples, \textsl{Designs, Codes and Cryptography}, 84 (2017), 87--100.

\bibitem{pless} V. Pless, A classification of self-orthogonal codes over $\text{GF}(2)$, \textsl{Discrete Mathematics},  3 (1972), 209--246.

\bibitem{sloanepless} V. Pless, N.\ J.\ A.\ Sloane, On the classification and enumeration of self-dual codes, \textsl{Journal of Combinatorial Theory (A)},  18 (1975), 313--335.

\bibitem{sloane} F.\ J.\ MacWilliams, N.\ J.\ A.\ Sloane, \textsl{The Theory of Error-Correcting Codes}, North-Holland (1983).


\bibitem{nakata} M.\ Nakata, A numerical evaluation of highly accurate multiple-precision arithmetic version of semidefinite programming solver: SDPA-GMP, -QD and -DD, \textsl{the proceedings of 2010 IEEE Multi-Conference on Systems and Control}, 2010, 29--34.


\bibitem{dreadnaut} B.\ D.\ McKay, A.\ Piperno, Practical graph isomorphism, II, \textsl{J.\ Symbolic Computation}, No.\ 60 (2013), 94--112.

\bibitem{cw4} S.\ C.\ Polak, Semidefinite programming bounds for constant weight codes, \textsl{IEEE Transactions on Information Theory} (2018), \href{https://doi.org/10.1109/TIT.2018.2854800}{\tt doi:10.1109/TIT.2018.2854800}

\bibitem{programs}  A.\ E.\ Brouwer, S.\ C.\ Polak, Computer programs, \url{https://github.com/codeuniqueness/code_uniqueness}. 


\bibitem{schrijver} A.\ Schrijver, New code upper bounds from the Terwilliger algebra and semidefinite programming, \textsl{IEEE Transactions on Information Theory}, 51 (2005), 2859--2866.


\bibitem{sdpa} M.\ Yamashita, K.\ Fujisawa, M.\ Fukuda, K.\ Kobayashi, K.\ Nakta, M.\ Nakata, Latest developments in the SDPA Family for solving large-scale SDPs, In \textsl{Handbook on Semidefinite, Cone and Polynomial Optimization: Theory, Algorithms, Software and Applications} edited by Miguel F. Anjos and Jean B. Lasserre, Springer, Chapter 24 (2011), 687--714. 

\end{thebibliography}
\end{document}